\documentclass[a4paper,11pt]{amsart}

  \def\cS{{\cal {S}}}
  
  \def\cN{{\cal {N}}}
\def\eps{{\varepsilon}}

\usepackage[dvips]{graphicx}
\usepackage{amssymb}
\usepackage{amsmath}
\usepackage[latin1]{inputenc}

\newtheorem{theo}{Theorem}[section]

\newtheorem{lem}[theo]{Lemma}

\newtheorem{coro}[theo]{Corollary}

\theoremstyle{definition}

\theoremstyle{remark} 

\theoremstyle{definition}

\newcommand{\cal}{\mathcal}

\begin{document}

\title{On the Maillet--Baker continued fractions}
\author[B.~Adamczewski, Y. Bugeaud]{Boris Adamczewski, Yann Bugeaud}

\begin{abstract}
We use the Schmidt Subspace Theorem 
to establish the transcendence of a class of quasi-periodic 
continued fractions. 
This improves earlier works of Maillet and of A. Baker. 
We also improve an old result of Davenport and Roth on the rate
of increase of the denominators of
the convergents to any real algebraic number.
\end{abstract}
\maketitle
\section{Introduction}


A central question in Diophantine approximation 
is concerned with how algebraic numbers can be approximated by rationals. 
This problem is intimately connected 
with the behaviour of their continued fraction expansion.  
In particular, it is widely believed that the 
continued fraction expansion of any irrational algebraic number 
$\xi$ either is eventually periodic 
(and we know that this is the case if, and only if, $\xi$ is a 
quadratic irrational), or it contains arbitrarily large partial quotients. 
Apparently,  this problem was first considered by 
Khintchine in  \cite{Khintchine} 
(we also refer the reader to \cite{Allouche,Shallit, Waldschmidt} 
for surveys including a discussion on this subject). 
Some speculations about the randomness of the continued fraction 
expansion of algebraic numbers of degree at least 
three have later been made by Lang \cite{Lang}. However,  
one shall admit that our knowledge on this topic is up 
to now very limited.

A first step consists in providing explicit 
examples of transcendental continued fractions. 
The first result of this type goes back to  the 
pioneering work of 
Liouville \cite{Liouville}, who constructed transcendental 
real numbers with a very fast growing sequence of partial quotients. 
Subsequently, 
various authors used deeper transcendence criteria from 
Diophantine approximation to 
construct other classes of transcendental continued fractions. 
Of parti\-cular interest 
is the work of Maillet \cite{Maillet} 
(see also Section 34 of Perron \cite{Perron}), 
who was the first to give explicit examples of transcendental 
continued fractions with bounded 
partial quotients. 
His work has later been carried on by A. Baker \cite{Baker62,Baker64}. 

More precisely, Maillet proved that if ${\bf a}=(a_n)_{n\geq 0}$ 
is a non-eventually periodic sequence of positive integers, and if  
there are infinitely many positive integers $n$ such that 
$$
a_n=a_{n+1}=\ldots =a_{n+\lambda(n)-1}, 
$$
then the real number $\xi=[a_0; a_1, a_2, \cdots]$ is transcendental, 
as soon as  $\lambda(n)$ is larger 
than a certain function of the denominator of the $n$-th
convergent to $\xi$. 
Actually, the result of Maillet is 
more general and also includes the case of repetitions 
of blocks of consecutive 
partial quotients (see Section \ref{resavant}). 
His proof is based on a general form of the Liouville inequality 
which limits the approximation of algebraic numbers by quadratic irrationals. 
Indeed, under the previous assumption, the quadratic irrational real numbers 
$\xi_n$, defined as having the eventually 
periodic continued fraction expansion 
$[a_0;a_1,\cdots,a_{n-1},a_n,a_n,a_n,\ldots]$, provide  
infinitely many 
`too good' approximations to $\xi$.

Not surprisingly, the breakthrough made by Roth 
in his 1955 seminal paper \cite{Roth55} 
leads to an improvement of this result. 
Indeed, Baker \cite{Baker62} used in 1962 
the Roth theorem for number fields obtained by LeVeque \cite{Leveque56} 
to strongly improve upon the results of 
Maillet and make them more explicit. 
His main idea was to remark that when infinitely many 
of the quadratic approximations found by Maillet lie 
in a same quadratic number field, one can favourably replace the use 
of the Liouville inequality by the one of LeVeque's Theorem. 

The purpose of the present paper is to improve 
the results obtained by Baker in \cite{Baker62}, that are recalled 
in Section \ref{resavant}.
Our approach rests on the Schmidt Subspace Theorem, 
but we follow a rather different way than the one previously 
considered by Maillet and by Baker. Our results are stated
in Section \ref{res} and proved in Section \ref{pr}. 
Section \ref{DR} is devoted to the improvement of
an old result of Davenport and Roth \cite{DavRoth} on the rate
of increase of the denominators of
the convergents to any real algebraic number. It is the key point
for the proof of Theorem \ref{amel4} below, and is also of
independent interest. Auxiliary results are gathered
in Section \ref{aux}.

\section{Earlier results}\label{resavant}

Throughout the present paper, we keep the following notation.
Let ${\bf a}=(a_n)_{n\geq 0}$ be a sequence of 
positive integers, that is not eventually periodic.
Let $(n_k)_{k\geq 0}$ be an increasing 
sequence of positive integers. Let
$(\lambda_k)_{k\geq 0}$ and $(r_k)_{k\geq 0}$ be sequences 
of positive integers. 
Assume that for any non-negative integer $k$, we have 
$n_{k+1} \ge n_k + \lambda_k r_k$ and
\begin{equation}\label{rep}
a_{m+r_k}=a_m\;\mbox{ for }\;  n_k\leq m\leq n_k+(\lambda_k - 1) r_k - 1,
\end{equation}
and consider the real number $\xi$ defined by 
$$
\xi=[a_0;a_1,a_2,\ldots,a_n,\ldots].
$$
Then, $\xi$ has a quasi-periodic continued fraction expansion
in the following sense: for any 
positive integer $k$, a block of $r_k$ consecutive partial 
quotients is repeated $\lambda_k$ times, 
such a repetition occurring just after the $(n_k-1)$-th partial quotient.

In \cite{Baker62}, Baker established three theorems, which strongly 
improved the pioneering work of Maillet.
The first one is very general.

\begin{theo}[A. Baker]\label{baker0}
With the previous notation, let us assume that 
\begin{equation}\label{borne}
\limsup_{k\to\infty}
\frac{r_k}{n_k} < + \infty.
\end{equation}
and
\begin{equation}\label{hypgen}
\limsup_{k\to\infty}
\frac{(\log \lambda_k) (\log n_k)^{1/2}}{n_k} = + \infty.
\end{equation}
Then, the real number $\xi$ is transcendental.
\end{theo}

Actually, it is not difficult 
to modify Baker's proof of Theorem~\ref{baker0}
in order to get rid of the assumption (\ref{borne}).
Notice also a related result due to Mignotte \cite{Mign}.
Under the assumption that the sequence ${\bf a}$ is bounded,
condition~(\ref{hypgen}) can be considerably relaxed.

\begin{theo}[A. Baker]\label{baker1}
Let $A \ge 2$ be an integer. Let ${\bf a}$ be a sequence of integers
at most equal to $A$ that satisfy (\ref{rep}) for a
bounded sequence $(r_k)_{k \ge 0}$. 
Assume that 
$$
\limsup_{k\to\infty}\frac{\lambda_k}{n_k}>B=B(A),
$$
where $B$ is defined by 
$$
B=2\left(\frac{\log\left( \left(A+\sqrt{A^2+4}\right)/2\right)}
{\log\left( (1+\sqrt{5})/2\right)}\right)-1.
$$
Then, the real number $\xi$ is transcendental.
\end{theo}

First, we remark that $B(A)$ increases with $A$ 
and that $\lim_{A\to\infty}B(A)=+\infty$. 
The smallest value, obtained for $A=2$, is $B(2)\simeq 2.66...$ 
Let us also note that, 
when one only knows that the sequence ${\bf a }$ 
is bounded, but without having any explicit 
bound, the stronger assumption 
\begin{equation}\label{infty}
\limsup_{k\to\infty}\frac{\lambda_k}{n_k}=+\infty
\end{equation}
is required  to apply Theorem \ref{baker1}.  

\medskip

One of the difficulties in the proof of Theorem \ref{baker1} 
is that one needs a precise estimate for 
the growth of the sequence of
the denominators of the convergents to $\xi$.
This in particular explains why, in this result, 
the value of $B$ depends on $A$.  
However, for a more restricted class of quasi-periodic continued fractions, 
that we present now, Baker \cite{Baker62} partly succeeded 
in overcoming this difficulty.

\begin{theo}[A. Baker]\label{baker2}
Let us consider the quasi-periodic continued fraction 
$$
\xi=[a_0;a_1,\ldots,a_{n_0-1},
\underbrace{a_{n_0},\ldots,a_{n_0+r_0-1}}_{\lambda_0}
,\underbrace{a_{n_1},\ldots,a_{n_1+r_1-1}}_{\lambda_1},\ldots],
$$
where the notation implies that $n_{k+1}=n_k+\lambda_kr_k$ and 
the $\lambda$'s indicate 
the number of times a block of partial quotients is repeated. 
Let us assume that the sequences 
$(a_n)_{n\geq 0}$ and $(r_k)_{k\geq 0}$ are both bounded, that
$(a_n)_{n \ge 0}$ is not ultimately periodic, and that 
\begin{equation}\label{ba}
\liminf_{k\to\infty}\frac{\lambda_{k+1}}{\lambda_k}>2.
\end{equation}
Then, the real number $\xi$ is transcendental.
\end{theo}

As a typical example of such continued fractions, 
Baker considered at the end of \cite{Baker62} 
the following family of real numbers: 
\begin{equation}\label{exba}
\xi_{a,b}=[0;\underbrace{a,a,\ldots,a}_{\lambda_0},
\underbrace{b,b,\ldots,b}_{\lambda_1},\underbrace{a,a,\ldots,a}
_{\lambda_2},b,b,\ldots],
\end{equation}
where $a$ and $b$ denote distinct positive integers.  
In the very particular case where $a=1$ and $b=2$, 
Baker improved Theorem \ref{baker2} by showing that 
$\xi$ is a transcendental number as soon as 
$\liminf_{k\to\infty}(\lambda_{k+1}/\lambda_k)>1.72$. 
Unfortunately, Baker's approach does not enable us to replace $2$ 
by a constant smaller than $\sqrt 2\simeq 1,41$ 
in Inequality (\ref{ba}), even for the  specific examples 
considered in (\ref{exba}).

\section{Main results}\label{res}

We present here our main results which 
improve the three theorems due to Baker 
mentioned in the previous Section.

\medskip

The first of Baker's results, namely Theorem \ref{baker0},  
heavily rests on an upper bound due to
Davenport and Roth \cite{DavRoth} (see (\ref{LiouvDR}) below) for the  
rate of increase of the denominators of
the convergents to any real algebraic number. Our improvement
of (\ref{LiouvDR}) stated in Theorem \ref{amelDR} below
allows us to get the following
strengthening of Theorem \ref{baker0}.

\begin{theo}\label{amel4}
Let ${\bf a}=(a_n)_{n\geq 0}$ be a sequence of 
positive integers, which satisfies (\ref{rep})  
and is not ultimately periodic. 
Assume that  
\begin{equation}\label{hypgenamel}
\limsup_{k\to\infty}
\frac{\log \lambda_k}{n_k^{\eps + 2/3}} = + \infty 
\end{equation}
holds for some $\eps > 0$. 
Then, the real number $\xi=[a_0;a_1,a_2,\ldots,a_n,\ldots]$ 
is transcendental.
\end{theo}

\medskip

In order to improve the two other results quoted in Section \ref{resavant}, 
it is tempting to try to apply the powerful Schmidt Subspace 
Theorem (see Section \ref{aux}) instead of the result 
of LeVeque mentioned in the Introduction.
For instance, the authors of \cite{Adamczewski_Bugeaud_Pal} 
recently improved Theorem~\ref{baker2} {\it via} the Subspace Theorem, 
but only in the particular case given in (\ref{exba}), for which 
they reached the bound $\sqrt 2$ (instead of $2$),
independently of the values of the distinct positive integers $a$ and $b$. 
See also related results by Davison \cite{Dav}.

\medskip

Quite surprisingly, a different application 
of the Subspace Theorem  based on the mirror 
formula (see Lemma \ref{lem4} for a definition) allows us to 
considerably relax the assumptions of 
two of the transcendence criteria obtained by Baker.
Our main result can be stated as follows.

\begin{theo}\label{amel1} 
Let ${\bf a}=(a_n)_{n\geq 0}$ be a sequence of 
positive integers, which satisfies (\ref{rep}) and is not ultimately periodic. 
Let $(p_n/q_n)_{n\geq 1}$ denote the sequence 
of convergents to the real number
$$
\xi=[a_0;a_1,a_2,\ldots,a_n,\ldots].
$$ 
Assume that the sequence $(q_n^{1/n})_{n\geq 1}$ is bounded 
(which is in particular the case when the sequence ${\bf a}$ is bounded), 
and that  
\begin{equation}\label{hyp}
\limsup_{k\to\infty}\frac{\lambda_k}{n_k}>0.
\end{equation} 
Then, the real number $\xi$ is transcendental.
\end{theo} 

Unlike in Theorem \ref{baker1}, 
the transcendence condition obtained in Theorem \ref{amel1} does not 
require neither that the partial quotients of the real number $\xi$
are bounded, nor that the lengths of the blocks which are repeated
are bounded. Furthermore, we point out that the assumption
`the sequence $(q_n^{1/n})_{n\geq 1}$ is bounded' is satisfied by
almost all real numbers. If one follows 
Baker's proof of Theorem \ref{baker0}
under this additional assumption, it is easily seen that one
gets a much weaker version of Theorem \ref{amel1}, namely with the  
condition (\ref{hyp}) being replaced by (\ref{infty}).

The proof of Theorem \ref{amel1} splits into two parts.
In the first part, we develop a new application of 
the Schmidt Subspace Theorem, 
based on Lemma \ref{lem4} below. This is the main novelty of
the present paper and it allows us to deal e.g.
with real numbers $\xi$ satisfying the assumption
of Baker's Theorem \ref{baker1}. 
The second part is far much easier.

\medskip

As a direct corollary of Theorem \ref{amel1}, 
we obtain the following improvement of Theorem \ref{baker2}.

\begin{coro}\label{amel2}
Let us consider the quasi-periodic continued fraction 
$$
\xi=[a_0;a_1,\ldots,a_{n_0-1},
\underbrace{a_{n_0},\ldots,a_{n_0+r_0-1}}_{\lambda_0}
,\underbrace{a_{n_1},\ldots,a_{n_1+r_1-1}}_{\lambda_1},\ldots],
$$
where the notation implies that $n_{k+1}=n_k+\lambda_kr_k$ and the 
$\lambda$'s indicate 
the number of times a block of partial quotients is repeated. 
Denote by $(p_n / q_n)_{n \ge 0}$ the sequence of the convergents to $\xi$.
Assume that the sequences $(q_n^{1/n})_{n \ge 0}$ and
$(r_k)_{k\geq 0}$ are bounded, that $(a_n)_{n \ge 0}$
is not ultimately periodic, and that
\begin{equation}\label{hyp2}
\liminf_{k\to\infty}\frac{\lambda_{k+1}}{\lambda_k}>1.
\end{equation} 
Then, the real number $\xi$ is transcendental.
\end{coro} 

Finally, we mention that 
applying the Schmidt Subspace Theorem 
in a similar way as in
our previous work \cite{Adamczewski_Bugeaud_Rep}
allows us to get rid of the assumptions on the sequences
$(a_n)_{n \ge 0}$ and $(r_k)_{k \ge 0}$ in Theorem \ref{baker2}.

\begin{theo}\label{amel3}
Let us consider the quasi-periodic continued fraction 
$$
\xi=[a_0;a_1,\ldots,a_{n_0-1},
\underbrace{a_{n_0},\ldots,a_{n_0+r_0-1}}_{\lambda_0}
,\underbrace{a_{n_1},\ldots,a_{n_1+r_1-1}}_{\lambda_1},\ldots].
$$
Assume that the sequence 
$(a_n)_{n\geq 0}$ is not ultimately periodic, and that 
\begin{equation}\label{bab}
\liminf_{k\to\infty}\frac{\lambda_{k+1}}{\lambda_k}>2.
\end{equation}
Then, the real number $\xi$ is transcendental.
\end{theo}


\section{An improvement of a result of Davenport and Roth}\label{DR}

Throughout the present Section (which can be read
independently of the rest of the paper),
$\xi$ denotes an arbitrary irrational, 
real algebraic number and $(p_n/q_n)_{n \ge 1}$ always denotes
the sequence of its convergents. The rate 
of growth of $(q_n)_{n \ge 1}$ is at least
exponential, as immediately follows from the theory
of continued fraction, see Lemma \ref{lem3} below.  
Our purpose is to estimate
it from above. It is well known that,
if $\xi$ is quadratic, then there
exists a real number $C(\xi)$, depending only on $\xi$, 
such that $q_n^{1/n} \le C(\xi)$ for any $n \ge 1$.
It is widely believed that $(q_n^{1/n})_{n \ge 1}$ also remains
bounded if the degree of $\xi$ is greater than two. However, we
seem to be very far away from a proof (or a disproof).

The first general upper estimate for the rate of increase
of $(q_n)_{n \ge 1}$ follows from the Liouville inequality,
saying that any algebraic number of degree $d$ cannot be
approximated by rationals at an order greater than $d$. Using this result,
we easily get that
\begin{equation}\label{Liouv}
\log \log q_n \ll n.  
\end{equation}
Throughout the present Section, all
the constants implied by $\ll$ depend only on $\xi$.

Let $\delta$ be a positive real number. In 1955, 
Roth \cite{Roth55} proved that the
set of solutions to the inequality
$$
\biggl| \xi - \frac{p}{q} \biggr| < \frac{1}{q^{2 + \delta}},
$$
in integers $p$, $q$ with $\gcd (p, q) = 1$ and $q > 0$, is finite.
In his joint work with Davenport \cite{DavRoth}, 
some steps from \cite{Roth55} were made more explicit  
in order to get an upper estimate for the
cardinality $\cN(\xi, \delta)$ of this set.
In particular, Davenport and Roth \cite{DavRoth} established that,
for $\delta \le 1/3$,
there exist positive constants $c_1$ and $c_2$, depending only
on $\xi$, such that 
\begin{equation}\label{DRoth}
\cN(\xi, \delta) \le c_1 \, \exp\{c_2 \delta^{-2} \}. 
\end{equation}
They further derived from (\ref{DRoth}) an
improvement of (\ref{Liouv}), namely the upper estimate
\begin{equation}\label{LiouvDR}
\log \log q_n \ll \frac{n }{\sqrt{\log n}}.  
\end{equation}

Bombieri and van der Poorten \cite{BombierivdP} 
were the first who established 
an upper bound for
$\cN(\xi, \delta)$ which is polynomial in $\delta^{-1}$.
A slight sharpening has subsequently
been obtained by Evertse, who proved at the end of 
Section 6 of \cite{Evertse97}
that, for $\delta < 1$, 
there exists a positive constant $c_3$, depending only
on $\xi$, such that 
\begin{equation}\label{Evert}
\cN(\xi, \delta) \le c_3 \, \delta^{-3} \, (1 + \log \delta^{-1})^2. 
\end{equation}
Any qualitative improvement of (\ref{DRoth}) yields an
improvement of (\ref{LiouvDR}). In particular, 
if we insert (\ref{Evert}) instead of (\ref{DRoth})
in Davenport and Roth's proof of (\ref{LiouvDR}), we get the upper estimate
\begin{equation}\label{LiouvEv}
\log \log q_n \ll n^{3/4} \, \sqrt{\log n}.  
\end{equation}
It turns out that a suitable modification of the argument used by
Davenport and Roth allows us to derive from (\ref{Evert}) a much better
result than (\ref{LiouvEv}).

\begin{theo}\label{amelDR}
Let $\xi$ be an arbitrary irrational, 
real algebraic number and let $(p_n/q_n)_{n \ge 1}$ denote
the sequence of its convergents.
Then, for any $\eps > 0$, there exists a constant $c_4$,
depending only on $\xi$ and $\eps$, such that
$$
\log \log q_n \le c_4 \,  n^{2/3 + \eps}.  
$$
\end{theo}

As an immediate corollary, we get a transcendence criterion for 
real numbers whose convergents have very large denominators.

\begin{coro}\label{cortrans}
Let $\theta$ be an irrational, 
real number and let $(r_n/s_n)_{n \ge 1}$ denote
the sequence of its convergents. If there exists
a positive real number $\eps$ such that
$$
\limsup_{n \to + \infty} \, \frac{\log \log s_n}{  n^{2/3 + \eps}} = + \infty,
$$
then $\theta$ is transcendental.
\end{coro}

Corollary \ref{cortrans} is the key point for
the proof of Theorem \ref{amel4}.

\medskip

\begin{proof}[Proof of Theorem \ref{amelDR}]
The basic idea is to introduce
more parameters in the proof of Theorem 3 of \cite{DavRoth}.
Recall that we have
\begin{equation}\label{frcont}
\biggl| \xi - \frac{p_n }{ q_n} \biggr| < \frac{1 }{ q_n q_{n+1}},  
\end{equation}
for any $n \ge 1$. 
Let $k \ge 1$ be an integer and $\delta_1, \ldots, \delta_k$ be real
numbers with $0 < \delta_1 < \delta_2 < \ldots < \delta_k < 1$, that will be
selected later on.

It is convenient to introduce a positive real number $\nu > 1$ such that 
$\cN(\xi, \delta) \ll \delta^{-\nu}$ holds for any 
$\delta$ with $0 < \delta < 1$.
In view of (\ref{Evert}), we can take for $\nu$ any real number
strictly larger than $3$.

Let $N$ be a (sufficiently large) integer
and put $\cS_0 = \{1, 2, \ldots, N\}$.
For $j = 1, \ldots, k$, let $\cS_j$ denote the set of positive
integers $n$ in $\cS_0$ 
such that $q_{n+1} > q_n^{1+ \delta_j}$.  Observe that
$\cS_0 \supset \cS_1 \supset \ldots \supset \cS_k$.
It follows from (\ref{frcont}) that, for any $n$
in $\cS_j$, the convergent $p_n / q_n$ 
gives a solution to
$$
\biggl| \xi - \frac{p}{  q} \biggr| < \frac{1}{  q^{2 + \delta_j}}.
$$
Consequently, the cardinality of $\cS_j$ is at most
$\cN(\xi, \delta_j)$, thus it is $\ll \delta_j^{-\nu}$.

Write
$$
\cS_0 = (\cS_0 \setminus \cS_1) \cup (\cS_1 \setminus \cS_2)  \cup \ldots
\cup (\cS_{k-1} \setminus \cS_k) \cup \cS_k.
$$
Let $j$ be an integer with $1 \le j \le k$.
The cardinality of $\cS_0 \setminus \cS_1$ is obviously bounded by $N$ and,
if $j \ge 2$, the cardinality of
$\cS_{j-1} \setminus \cS_j$ is $\ll \delta_{j-1}^{-\nu}$.
Furthermore, for any $n$ in 
$\cS_{j-1} \setminus \cS_j$, we get
$$
\frac{\log q_{n+1}}{  \log q_n} \le 1 + \delta_j.
$$
Denoting by $d$ the degree of $\xi$, we infer from (\ref{frcont}) and the
Liouville inequality that
$$
\frac{\log q_{n+1}}{  \log q_n} \le d
$$
holds for every sufficiently large
integer $n$ in $\cS_k$. Combining these estimates
with the fact that $\cS_k$ has $\ll \delta_k^{-\nu}$ elements,
we obtain that
$$
\log q_N \ll \frac{\log q_N }{ \log q_{N-1}} \times
\frac{\log q_{N-1}}{  \log q_{N-2}}  \times \ldots \times 
\frac{\log q_3}{ \log q_2}  
\ll (1 + \delta_1)^N \, \prod_{j=2}^{k} \, 
(1 + \delta_{j})^{\delta_{j-1}^{-\nu}} 
\, d^{\delta_k^{-\nu}}. 
$$
Taking the logarithm, we get
\begin{equation}\label{bornelog}
\log \log q_N \ll N \log (1 + \delta_1) + \sum_{j = 2}^{k} \,
\delta_{j-1}^{-\nu} \log (1 + \delta_{j}) + \delta_k^{-\nu}.
\end{equation}
We now select $\delta_1, \ldots, \delta_k$. For $j = 1, \ldots, k$, set
$$
\delta_j = N^{- (\nu^k - \nu^{j-1})/(\nu^{k+1} - 1)}.
$$
We check that $0 < \delta_1 < \ldots < \delta_k < 1$, and we 
easily infer from (\ref{bornelog}) that
\begin{equation}\label{bornenu}
\log \log q_N \ll k \, N^{(\nu^{k+1} - \nu^k)/(\nu^{k+1} - 1)}
= k \, N^{(\nu - 1)/(\nu - \nu^{-k})}.  
\end{equation}
By (\ref{Liouv}), we may assume that $\eps \le 1/3$.
In view of (\ref{Evert}), we can take $\nu = 3/(1 - \eps)$.
Choosing then for $k$ the smallest 
integer greater than $\log \eps^{-1}$,
we get from (\ref{bornenu}) that
$$
\log \log q_N \ll (\log \eps^{-1}) \, N^{\eps + 2/3},
$$
as claimed. 
\end{proof}


\section{Auxiliary results}\label{aux}

Our Theorems \ref{amel1} and \ref{amel3} rest on the 
Schmidt Subspace Theorem \cite{Schmidt72a} 
(see also \cite{Schmidt80}), that we recall now. 

\begin{theo}[W. M. Schmidt]\label{sousespace} 
Let $m \ge 2$ be an integer.
Let $L_1, \ldots, L_m$ be linearly independent linear forms in
${\bf x} = (x_1, \ldots, x_m)$ with algebraic coefficients.
Let $\varepsilon$ be a positive real number.
Then, the set of solutions ${\bf x} = (x_1, \ldots, x_m)$ 
in ${\mathbb Z}^m$ to the inequality
$$
\vert L_1 ({\bf x}) \ldots L_m ({\bf x}) \vert  \le
(\max\{|x_1|, \ldots , |x_m|\})^{-\varepsilon}
$$
lies in finitely many proper subspaces of ${\mathbb Q}^m$.
\end{theo}

\medskip

For the reader convenience, we recall here some classical 
results from the theory of continued 
fractions, whose proofs can be found for example 
in the book of Perron \cite{Perron}.

\medskip

\begin{lem}\label{lem1}
Let $\xi = [a_0; a_1, a_2, \cdots]$ and 
$\eta = [b_0; b_1, b_2, \cdots]$ be real numbers. 
Let $n \geq 1$ such that $a_j = b_j$ for any $j=0,\ldots, n$. We then have
$\vert\xi - \eta\vert \le q_n^{-2}$, where $q_n$ denotes the denominator
of the $n$-th convergent to $\xi$.
\end{lem}

\medskip

\begin{lem}\label{lem3}
Let $(a_n)_{n \ge 0}$ be a sequence of 
positive integers at most equal to $M\!$, 
let $n$ be a positive integer and set 
$p_n/q_n=[a_0;a_1,a_2,\ldots,a_n]$. Then, we have
$$
\sqrt{2}^{(n-1)} \le q_n \le (M+1)^n.
$$
\end{lem}

\medskip

The following innocent-looking formula appears to be
the key point in the proof of Theorem \ref{amel1}. 
In what follows, Equality (\ref{mir}) 
will be referred to as the mirror formula.

\begin{lem}\label{lem4}
Let $\xi = [a_0; a_1, a_2, \cdots]$ be 
a real number with convergents $(p_n/q_n)_{n\geq 1}$. 
Then, for any $n\geq 2$, we have 
\begin{equation}\label{mir}
\frac{q_n}{q_{n-1}}=[a_n;a_{n-1},\ldots,a_1].
\end{equation}
\end{lem}

\medskip

For positive integers $a_1, \ldots, a_m$, denote 
by $K_m (a_1, \ldots, a_m)$ the denominator of the rational number
$[0; a_1, \ldots, a_m]$. It is commonly called a {\it continuant}.

\medskip

\begin{lem}\label{lem5}
For any positive integers $a_1, \ldots, a_m$ and any integer $k$ with
$1 \le k \le m-1$, we have
$$
K_m (a_1, \ldots , a_m) = K_m (a_m, \ldots, a_1),
$$

$$
K_k (a_1, \ldots, a_k) \cdot  K_{m-k} (a_{k+1}, \ldots, a_m)
\le K_m (a_1, \ldots , a_m)$$

$$\hspace{2.5cm}\le 2 \, K_k (a_1, \ldots, a_k) \cdot 
K_{m-k} (a_{k+1}, \ldots, a_m)
$$

\noindent and

$$\begin{array}{ll}
\frac{1}{2} \, K_m (a_k, \ldots, a_m, a_1, \ldots, a_{k-1}) &\le
K_m (a_1, \ldots, a_m) \\ \\
&\le 2 \, K_m (a_k, \ldots, a_m, a_1, \ldots, a_{k-1}).
\end{array}$$
\end{lem}

We finish our series of lemmas by an immediate consequence of
Roth's theorem.

\begin{lem}\label{lem6}
Let $(p_n / q_n)_{n \ge 0}$ denote the sequence of partial quotients
of a real number $\xi$. Let $\eta$ be a positive integer.
If $\xi$ is algebraic, then $q_{n+1} \le q_n^{1+ \eta}$
holds for any integer $n$ sufficiently large.
\end{lem}


\section{Proofs of our main results}\label{pr}


\begin{proof}[Proof of Theorem \ref{amel4}]
We follow the proof of Theorem 1 
from~\cite{Baker62}, except that we use Theorem \ref{amelDR}
instead of (\ref{LiouvDR}) and that we suitably apply Lemma \ref{lem5}
to get rid of the assumption (\ref{borne}). For completeness, we give
the details of the argument.

Assume that $\xi$ is algebraic of degree $d$.
For any positive integer $k$, set
$$
\xi_k := [a_0;a_1,\cdots,a_{n_k-1},
\underbrace{a_{n_k}, \ldots, a_{n_k + r_k - 1}}_{\infty}].
$$
Since the height of $\xi_k$ is at most $2 q_{n_k + r_k - 1}^2$,
the Liouville inequality 
(see e.g. \cite{BuLiv}, Corollary A.2) 
and Lemma \ref{lem5} give us that 
$$
| \xi - \xi_k| \gg q_{n_k + r_k - 1}^{-2d} \gg 
q_{n_k-1}^{-2d} \, K_{r_k} (a_{n_k}, \ldots, a_{n_k + r_k - 1})^{-2d}.
$$
Here and below, the constants implied by $\ll$ depend only on $\xi$. 
However, we infer from Lemmas \ref{lem1}, \ref{lem3} and \ref{lem5} that
$$
\begin{array}{ll}
 \hspace{-0.30cm} |\xi - \xi_k| 
& \hspace{-0.25cm}\ll q_{n_k + \lambda_k r_k - 1}^{-2} \\ \\
&\hspace{-0.4cm}\ \ll q_{n_k-1}^{-2} \, K_{r_k} (a_{n_k}, \ldots, a_{n_k + r_k - 1})^{-2d} \, 
K_{2 r_k} (a_{n_k}, \ldots, a_{n_k + 2 r_k - 1})^{- \lambda_k + 2d}
\\ \\
&\hspace{-0.4cm}\ \ll q_{n_k-1}^{-2} \, K_{r_k} (a_{n_k}, \ldots, a_{n_k + r_k - 1})^{-2d} \, 
2^{- \lambda_k / 2},
\end{array}
$$
if $k$ is sufficiently large.
A combination of the last two inequalities gives that
$$
\lambda_k \ll \log q_{n_k}.
$$
By the assumption (\ref{hypgenamel}), we then get that
$$
\limsup_{k\to\infty}
\frac{\log \log q_{n_k}}{n_k^{\eps + 2/3}} = + \infty,
$$
a contradiction with Corollary \ref{cortrans}.
\end{proof}


\begin{proof}[Proof of Theorem \ref{amel1}]
For any $k \ge 0$, set 
$$
K_k := K_{r_k} (a_{n_k}, \ldots , a_{n_k + r_k - 1}).
$$
By assumption, there exist a positive real number $c$
and an infinite set of integers ${\cal K}_1$,
ranged in increasing order, such that $\lambda_k \ge c n_k$ 
for any $k$ in ${\cal K}_1$.

The proof splits into two parts. Assume first
that the sequence $(K_k)_{k \in {\cal K}_1}$ is bounded. Since the $K_k$'s are 
non-negative integers, it follows that infinitely of them take the same value.   
Then, Lemma \ref{lem3} implies that there exist a positive integer $r$,
positive integers $b_0, \ldots, b_{r-1}$ and an infinite set ${\cal K}_2$
of positive integers such that
$$
r_k = r, \qquad a_{n_k + j} = b_j, \quad (0 \le j \le r-1),
$$
for any $k$ in ${\cal K}_2$.

Let $\alpha$ denote the real number having the 
purely periodic continued fraction expansion
with period $B = (b_{r-1}, \ldots , b_0)$, that is,
$$
\alpha=[b_{r-1} ; b_{r-2}, \ldots, b_0, b_{r-1}, 
\ldots , b_0, b_{r-1}, \ldots] = [B , B, \ldots, B, \ldots].
$$
Then, $\alpha$ is a quadratic number. 
We need to introduce some more notation. 
Let us denote by $p_n/q_n$ (respectively, by $r_n/s_n$) 
the $n$-th convergent to $\xi$ (respectively, to $\alpha$).  
Then, for any $k$ in ${\cal K}_2$,
set  $P_k = p_{n_k + \lambda_k r_k - 1},\; 
Q_k=q_{n_k + \lambda_k r_k - 1},\;
P'_k=p_{n_k + \lambda_k r_k - 2},\; Q'_k=q_{n_k + \lambda_k r_k - 2}$ and 
$S_k=s_{r\lambda_k-1}$. 

\medskip

By assumption, we already know that $\xi$ is irrational and not quadratic. 
Therefore, we assume that $\xi$ is algebraic and we aim at deriving a 
contradiction.

Let $k$ be in ${\cal K}_2$. 
By the theory of continued fractions, we have 
\begin{equation}\label{b1}
\vert Q_k\xi-P_k\vert<
\frac{1}{Q_k}\;\mbox{ and }\; 
\vert Q'_k\xi-P'_k\vert<\frac{1}
{Q'_k}.
\end{equation}
On the other hand, since by assumption 
$$
\frac{P_k}{Q_k}=[a_0;a_1,\cdots,a_{n_k-1},
\underbrace{B, B,\ldots,B}_{\lambda_k}],
$$
we get from the mirror formula (see Lemma \ref{lem4}) that 
$$
\frac{Q_k}{Q'_k}=[\underbrace{B,B,\ldots,B}_{\lambda_k},a_{n_k-1},\cdots,a_1].
$$
Then, Lemma~\ref{lem1} implies 
\begin{equation}\label{b2}
\left\vert Q'_k\alpha-Q_k\right\vert<\frac{Q'_k}{S_k^2}
\end{equation}
and we {\it a fortiori} obtain that 
\begin{equation}\label{bb}
\lim_{{\cal K}_2 \ni k\to\infty} \frac{Q_k}{Q'_k}=\alpha.
\end{equation}

Consider now the four linearly independent linear forms, 
whose coefficients are by assumption algebraic:
$$
\begin{array}{l}
L_1(X_1, X_2, X_3, X_4)=\xi X_1-X_3,\\
L_2(X_1, X_2, X_3, X_4)=\xi X_2 - X_4,\\ 
L_3(X_1, X_2, X_3, X_4)=\alpha X_2-X_1,\\
L_4(X_1, X_2, X_3, X_4)=X_1.
\end{array}
$$

Evaluating them on the quadruple $(Q_k,Q'_k,P_k,P'_k)$, it follows from 
(\ref{b1}) and (\ref{b2})  that
\begin{equation}\label{b3}
\prod_{1 \le j \le 4} \, |L_j (Q_k,Q'_k,P_k,P'_k)|
< \frac{1}{S_k^2}. 
\end{equation}
Let $M$ be an 
upper bound for the sequence $(q_n^{1/n})_{n\geq 1}$.
We infer from Lemma \ref{lem3} that $Q_k\leq (M+1)^{n_k+r\lambda_k}$ and 
$S_k\geq (\sqrt{2})^{r\lambda_k-2}$, 
for any positive integer $k$ in ${\cal K}_2$.
It thus follows that   
$$
S_k\geq (M+1)^{\left(\frac{\log{\sqrt 2}}{\log (M+1)}\right)(r\lambda_k-2)}
\geq 
Q_k^{\left(\frac{\log{\sqrt 2}}{\log (M+1)}\right)
\cdot\left(\frac{r\lambda_k-2}{n_k+r\lambda_k}\right)},
$$ 
for any positive integer $k$ in ${\cal K}_2$.
In particular, we get from (\ref{b3}) and (\ref{hyp}) that 
$$
\prod_{1 \le j \le 4} \, |L_j (Q_k,Q'_k,P_k,P'_k)|
\leq Q_k^{-\varepsilon}
$$
holds for some positive real number $\varepsilon$ and for $k$ large enough
in ${\cal K}_2$.

It then follows from Theorem~\ref{sousespace} 
that the points $(Q_k,Q'_k,P_k,P'_k)$ for $k$ in ${\cal K}_2$
lie in a finite number of proper subspaces of ${\mathbb Q}^4$. 
Thus, there exist a nonzero integer quadruple $(x_1,x_2,x_3,x_4)$ and
an infinite set of distinct positive integers 
${\cal K}_3 \subset {\cal K}_2$ such that
\begin{equation}\label{b6}
x_1 Q_k + x_2 Q'_k + x_3 P_k + x_4 P'_k = 0, 
\end{equation}
for any $k$ in ${\cal K}_3$.
Dividing (\ref{b6}) by $Q'_k$, we obtain
\begin{equation}\label{b7}
x_1\frac{Q_k}{Q'_k} + x_2 + x_3 \frac{P_k}{Q_k}\cdot\frac{Q_k}{Q'_k}
+ x_4 \frac{P'_k}{Q'_k}= 0. 
\end{equation}
By letting $k$ tend to infinity along ${\cal K}_3$ in (\ref{b7}), 
we derive from (\ref{bb}) that   
$$
x_1\alpha+x_2+(x_3\alpha+x_4)\xi=0.
$$
Since $\xi$ is not quadratic, it {\it a fortiori} cannot lie in 
${\mathbb Q}(\alpha)$. This implies that $x_3\alpha+x_4=0$ and, since 
$\alpha$ is irrational, it follows that $x_3=x_4=0$. Then,   
again by using that $\alpha$ is irrational, 
we get that $x_1=x_2=x_3=x_4=0$, 
which is a contradiction. This concludes the proof when the
sequence $(K_k)_{k \in {\cal K}_1}$ is bounded.

\medskip

Assume now that the sequence $(K_k)_{k \in {\cal K}_1}$ is unbounded.
Then, there exists an infinite set ${\cal K}_4$ of integers,
ranged in increasing order, such that the sequence
$(K_k)_{k \in {\cal K}_4}$ increases to infinity.

Recall that $(p_n / q_n)_{n \ge 0}$
denotes the sequence of convergents to $\xi$
and that $M$ denotes an upper bound for the
sequence $(q_n^{1/n})_{n \ge 0}$.
Let $d$ be a positive integer.
Let $k$ be in ${\cal K}_4$ and large enough in order 
that $\lambda_k \ge d + 1$ and 
\begin{equation}\label{condMc}
K_k \ge   M^{2d/c},
\end{equation}
with the constant $c$ defined at the beginning of the proof.
Then, the real number $\xi$ is very close
to the quadratic number
$$
\xi_k := [a_0;a_1,\cdots,a_{n_k-1},
\underbrace{a_{n_k}, \ldots, a_{n_k + r_k - 1}}_{\infty}].
$$
Define the polynomial
\begin{multline*}
 P_k (X):= (q_{n_k-2} q_{n_k + r_k-1} - 
q_{n_k-1} q_{n_k + r_k-2} ) X^2 \\ - (q_{n_k-2} p_{n_k + r_k-1} - q_{n_k-1} p_{n_k + r_k-2} 
+ p_{n_k-2} q_{n_k + r_k -1} - p_{n_k - 1} q_{n_k + r_k-2}) X  \\
+ (p_{n_k-2} p_{n_k + r_k -1} - p_{n_k -1} p_{n_k + r_k-2}),
\end{multline*}

\noindent and observe that $P_k (\xi_k) = 0$.
For any positive integer $k$, we infer from Rolle's Theorem 
and Lemma \ref{lem1} that
$$\begin{array}{ll}
|P_k (\xi)| &= |P_k (\xi) - P_k (\xi_k)|
\ll \, q_{n_k - 1} \,  q_{n_k + r_k -1} \, |\xi - \xi_k| \\ \\
& \ll 
q_{n_k - 1} \,  q_{n_k + r_k - 1}  \, q_{n_k + \lambda_k r_k -1}^{-2},
\end{array}$$
since the first $n_k + \lambda_k r_k - 1$ 
partial quotients of $\xi$ and $\xi_k$
are the same. Here and below, the constants implied
in $\ll$ depend at most on $\xi$ and on $d$, but they are
independent on $k$.
Now, it follows from Lemma \ref{lem5} that
$$
q_{n_k + \lambda_k r_k - 1} \ge q_{n_k - 1} \, K_k^{\lambda_k},
$$
thus, by (\ref{condMc}) and by Lemma \ref{lem5} again, we get
$$
|P_k (\xi)|  \ll K_k^{1 - 2 \lambda_k}
\ll K_k^{- \lambda_k - d} \ll (M^{2 n_k} \, K_k)^{-d} \ll 
(q_{n_k - 1} \,  q_{n_k + r_k -1})^{-d},
$$
since $\lambda_k \ge d+1$. Recalling
that $\xi$ is irrational and not quadratic,
it then follows
from the Liouville inequality (see e.g. \cite{BuLiv}, Theorem A.1) 
that $\xi$ cannot be algebraic of degree smaller than $d$.
Since $d$ is arbitrary,
this concludes the proof when the
sequence $(K_k)_{k \in {\cal K}_1}$ is bounded.
\end{proof}

\bigskip


\begin{proof}[Proof of Corollary \ref{amel2}]
Let us consider the quasi-periodic continued fraction 
$$
\xi=[a_0;a_1,\ldots,a_{n_0-1},
\underbrace{a_{n_0},\ldots,a_{n_0+r_0-1}}_{\lambda_0}
,\underbrace{a_{n_1},\ldots,a_{n_1+r_1-1}}_{\lambda_1},\ldots],
$$
satisfying the assumption of the corollary, and suppose that we have
\begin{equation}\label{h}
\liminf_{k\to\infty}\frac{\lambda_{k+1}}{\lambda_k}>1.
\end{equation}

For $k\geq 1$, we get that 
$$
n_k=n_0+ \sum_{j=0}^{k-1} r_j \lambda_j.
$$
Moreover, we infer from (\ref{h}) that there exist
positive real numbers $\delta$ and $M$ such that 
$$
\lambda_j<\frac{M \lambda_k}{(1+\delta)^{k-j}},
$$
for any $j<k$ with $k$ large enough. Since the sequence $(r_k)_{k\geq 0}$ 
is bounded, there exists a positive real number $r$ such that 
$$
n_k<n_0+r \lambda_k \, \sum_{j\geq 1}\frac{M}{(1+\delta)^j} 
$$
and thus 
$$
\limsup_{k \to \infty}
\frac{\lambda_k}{n_k}\geq \frac{1}{n_0+rM
\bigl(\sum_{j\geq 1}{(1+\delta)^{-j}}\bigr)}>0,$$
for $k$ large enough. Applying Theorem \ref{amel1}, this concludes the proof.
\end{proof}

\medskip


\begin{proof}[Proof of Theorem \ref{amel3}]
For any $k \ge 0$, set 
$$
K_k := K_{r_k} (a_{n_k}, \ldots , a_{n_k + r_k - 1}).
$$
In view of Corollary \ref{amel2}, there is no restriction
in assuming that the sequence $(K_k)_{k \ge 0}$ is unbounded.
Then, there exists an infinite set ${\cal K}_5$ of integers,
ranged in increasing order, such that the sequence
$(K_k)_{k \in {\cal K}_5}$ increases to infinity and such that,
for any $k$ in ${\cal K}_5$ and any integer $j$
with $0 \le j < k$, we have $K_j < K_k$.

Let $k$ be in ${\cal K}_5$.
The real number $\xi$ is very close
to the quadratic number
$$
\xi_k := [a_0;a_1,\cdots,a_{n_k-1},
\underbrace{a_{n_k}, \ldots, a_{n_k + r_k - 1}}_{\infty}].
$$
Let $j_k$ be the largest integer $< n_k$ such that
$a_{j_k} \not= a_{j_k + r_k}$. Choosing $k$ sufficiently
large, $j_k$ is well defined, since $(a_k)_{k \ge 0}$ is
not ultimately periodic. Observe that
$$
\xi_k := [a_0;a_1,\cdots,a_{j_k},
\underbrace{a_{j_k + 1}, \ldots, a_{j_k + r_k}}_{\infty}].
$$
Let $\eta \le 1$ be a positive real number.
We then assume that $\xi$ is algebraic and
we proceed as in the proof of Theorem 2 from
\cite{Adamczewski_Bugeaud_Rep}. 

Define the polynomial
$$\begin{array}{ll}
P_k (X) := & (q_{j_k-1} q_{j_k + r_k} - q_{j_k} q_{j_k + r_k-1} ) X^2  \\ \\
&- (q_{j_k-1} p_{j_k + r_k} - q_{j_k} p_{j_k + r_k-1} 
+ p_{j_k-1} q_{j_k + r_k} - p_{j_k} q_{j_k + r_k-1}) X \\ \\
& + (p_{j_k-1} p_{j_k + r_k} - p_{j_k} p_{j_k + r_k-1}),
\end{array}$$
and observe that $P_k (\xi_k) = 0$.
For any positive integer $k$ in ${\cal K}_5$, we infer from Rolle's Theorem 
and Lemma \ref{lem1} that
$$
|P_k (\xi)| = |P_k (\xi) - P_k (\xi_k)|
\ll \, q_{j_k} \,  q_{j_k + r_k} \, |\xi - \xi_k| \ll 
q_{j_k} \,  q_{j_k + r_k}  \, q_{j_k + \lambda_k r_k}^{-2},
$$
since the first $j_k + \lambda_k r_k$ partial quotients of $\xi$ and $\xi_k$
are the same.

Since $\xi$ is assumed to be algebraic, 
it follows from Lemma \ref{lem6} that
$$
|(q_{j_k-1} q_{j_k + r_k} - q_{j_k} q_{j_k + r_k-1} ) \xi
- (q_{j_k-1} p_{j_k + r_k} - q_{j_k} p_{j_k + r_k-1} )| \ll
q_{j_k} \, q_{j_k + r_k}^{-1 + \eta}  
$$
and
$$
|(q_{j_k-1} q_{j_k + r_k} - q_{j_k} q_{j_k + r_k-1}  ) \xi
- (p_{j_k-1} q_{j_k + r_k} - p_{j_k} q_{j_k + r_k-1})| \ll
q_{j_k}^{-1 + \eta} \, q_{j_k + r_k}, 
$$
if $k$ in ${\cal K}_5$ is large enough.
Furthermore, we have as well the obvious upper bound
$$
|q_{j_k-1} q_{j_k + r_k} - q_{j_k} q_{j_k + r_k -1} | 
\le  q_{j_k} \, q_{j_k + r_k}.
$$

Consider now the four linearly independent linear forms with algebraic 
coefficients:
$$
\begin{array}{l}
L_1(X_1, X_2, X_3, X_4) = \xi^2 X_1 - \xi (X_2 + X_3) + X_4,  \\
L_2(X_1, X_2, X_3, X_4) = \xi X_1 - X_2, \\
L_3(X_1, X_2, X_3, X_4) = \xi X_1 - X_3, \\
L_4(X_1, X_2, X_3, X_4) = X_1. 
\end{array}
$$
Evaluating them on the quadruple 
$$
{\underline z_k} := (q_{j_k-1} q_{j_k + r_k}  - q_{j_k} q_{j_k + r_k-1}, 
q_{j_k-1} p_{j_k + r_k} - q_{j_k} p_{j_k + r_k-1},  
$$
$$
p_{j_k-1} q_{j_k + r_k} - p_{j_k} q_{j_k + r_k-1}, 
p_{j_k-1} p_{j_k + r_k} - p_{j_k} p_{j_k + r_k-1}), 
$$ 
we find that
$$
\Pi := \prod_{1 \le j \le 4} \, |L_j ({\underline z_k})|
\ll (q_{j_k} \,  q_{j_k + r_k})^{2 + \eta}  \, q_{j_k + \lambda_k r_k}^{-2}.
$$
Now, it follows from the last assertion of Lemma \ref{lem5} that
$$
q_{j_k + \lambda_k r_k} \ge q_{j_k} \, (K_k / 2)^{\lambda_k},
$$
thus, by Lemma \ref{lem5} again,
$$
\Pi \ll q_{j_k}^{2 + 2 \eta} \, (2 K_k)^{2 + \eta} 
\, (K_k / 2)^{- 2 \lambda_k}
\ll K_k^{2 (1 + \eta)(\lambda_1 + \ldots + \lambda_{k-1} + k)}
\, (K_k / 2)^{- 2 \lambda_k}.
$$
Using hypotheses (\ref{bab}) and choosing $\eta$ small enough, we infer 
from the preceding inequality that there exists a positive real
number $\varepsilon$ such that
$$
\Pi \ll K_k^{- 2(\lambda_1 + \ldots + \lambda_{k-1} + k) \varepsilon} 
\ll (q_{j_k} \,  q_{j_k + r_k})^{- \varepsilon}.
$$
It then follows from Theorem~\ref{sousespace} 
that the points ${\underline z_k}$ for $k$ in ${\cal K}_5$
lie in a finite number of proper subspaces of ${\mathbb Q}^4$. 
Thus, there exist a nonzero integer quadruple $(x_1,x_2,x_3,x_4)$ and
an infinite set of distinct positive integers 
${\cal K}_6 \subset {\cal K}_5$ such that
$$
x_1 (q_{j_k-1} q_{j_k + r_k}  - q_{j_k} q_{j_k + r_k-1}) + 
x_2 (q_{j_k-1} p_{j_k + r_k} - q_{j_k} p_{j_k + r_k-1})
$$
$$
+ x_3 (p_{j_k-1} q_{j_k + r_k} - p_{j_k} q_{j_k + r_k-1})
+ x_4 (p_{j_k-1} p_{j_k + r_k} - p_{j_k} p_{j_k + r_k-1}) = 0. 
$$
for any $k$ in ${\cal K}_6$.
We then argue exactly as in the proof of Theorem 2 from
\cite{Adamczewski_Bugeaud_Rep}. This is made possible
by our choice of $j_k$. We then reach a contradiction,
which concludes the proof of our theorem.
\end{proof}


\vspace{1 cm}

\noindent Boris Adamczewski    \hfill{Yann Bugeaud}

\noindent  CNRS, Institut Camille Jordan 
\hfill{Universit\'e Louis Pasteur}

\noindent   Universit\'e Claude Bernard Lyon 1 
\hfill{U. F. R. de math\'ematiques}

\noindent   B\^at. Braconnier, 21 avenue Claude Bernard
 \hfill{7, rue Ren\'e Descartes}

\noindent   69622 VILLEURBANNE Cedex    
\hfill{67084 STRASBOURG Cedex}

\noindent FRANCE  \hfill{FRANCE}

\vskip1mm

\noindent {\tt Boris.Adamczewski@math.univ-lyon1.fr}
\hfill{{\tt bugeaud@math.u-strasbg.fr}}

\end{document}